# Letter counting: a stem cell for Cryptology, Quantitative Linguistics, and Statistics

Bernard Ycart
Univ. Grenoble Alpes
Laboratoire Jean Kuntzmann CNRS UMR 5224
51 rue des Mathématiques, 38041 Grenoble cedex 9, France
Tel: +33 4 76 51 49 95
Email: Bernard.Ycart@imag.fr

**Argument**

Counting letters in written texts is a very ancient practice. It has accompanied the development of Cryptology, Quantitative Linguistics, and Statistics. In Cryptology, counting frequencies of the different characters in an encrypted message is the basis of the so called frequency analysis method. In Quantitative Linguistics, the proportion of vowels to consonants in different languages was studied long before authorship attribution. In Statistics, the alternation vowel-consonants was the only example that Markov ever gave of his theory of chained events. A short history of letter counting is presented. The three domains, Cryptology, Quantitative Linguistics, and Statistics, are then examined, focusing on the interactions with the other two fields through letter counting. As a conclusion, the eclectism of past centuries scholars, their background in humanities, and their familiarity with cryptograms, are identified as contributing factors to the mutual enrichment process which is described here.

**Introduction**

In his founding "Manuscript on Deciphering Cryptographic Messages", Yaʿkūb ibn Isḥāq Al Kindī (ca. 801–873) wrote the following (Mrayati et al. 2002, 122):

> The quantitative expedients include determining the most frequently occurring letters in the language in which cryptograms are to be cryptanalysed. If vowels functioned as the material from which any language is made, and non-vowels functioned as the shape of any language, and since many shapes can be made from the same material, then the number of vowels in any language would be greater than non-vowels. For instance, gold is the material of many shapes of finery and vessels; it may cover crowns, bangles, cups, etc. The gold in these realizations is more than the shape made of it. Similarly, the vowels which are the material of any kind of text are more than the non-vowels in any language.

Six centuries later, Leon Battista Alberti (1404–1472) echoed (Alberti 2010, 173):

> I could go on at great length about this, but as briefly as possible, just let me say that the arrangement and the number of letters constitute various syllables which, when combined, form finished words, which differ in sound and meaning. First I shall speak of number and how this is connected to numeric ratios. This phenomenon above all involves the vowels. We will begin with them. A syllable is constituted of either a single vowel, a consonant associated to a vowel, or several consonants united to a vowel. Without vowels there are no syllables. Thus if we take, for example, one or two pages of poetry or prose and extract the vowels and consonants, listing them in separate series, vowels on one side and consonants on the other, you will no doubt find that there are numerous vowels.

There is no evidence that Alberti had any knowledge of earlier Arab cryptology treatises, which were rediscovered only recently (Al-Kadit 1992). Actually, many similar observations can be found in Cryptology treatises along the centuries. All have in common a very prosaic, long time occupation of scholars: counting letters. As Georges Udny Yule (1871–1951) said: "It is full of interest, and is of the most venerable antiquity and august associations[…]" (Yule 1944, 7).

Since Al Kindī, all Cryptology treatises have included a section on monoalphabetic simple substitution cyphers, which consist in replacing each letter by a given character, through a one-to-one correspondence. If the language of the encrypted message, and the frequency order of the letters in that language are known, the code is vulnerable to so called *frequency analysis*: the most frequent character in the message stands for the most frequent letter of the alphabet, the next more frequent for the second most frequent letter, and so on. This seemingly straightforward procedure conceals a major difficulty, which took a very long time to understand: there is no such thing as *the* frequency of letters. Along the centuries, many scholars have struggled with that difficulty, gradually developing notions that went far beyond technical necessities of cryptology, and eventually helped shaping two other fields: Quantitative Linguistics and Statistics. Throughout their history, the three disciplines have constantly nourished each other, and letter counting has been the catalyst of that mutual enrichment process: its history is the object of the present article.

The historical development of each of the three disciplines is well documented. For Cryptology, basic references are (Kahn 1996) and (Mollin 2005). Here, our focus will be on the perception of letter frequencies among the authors of Cryptology treatises, from Al Kindī to 20[th] century textbooks. Arab authors often gave tables of letter frequencies, based on counts of several thousand letters, and they appeared quite aware of the necessity of a large sample size. Western treatises after Alberti contrast: quantified indications remain scarce, and large scale counts do not appear before the 19[th] century. Alberti's "De Componendis Cifris" is a quite singular point in that history: it went notably further than cryptology, giving the first example of a quantified stylistic difference ever. After Alberti, the emergence of Quantitative Linguistics was slow. Before the 19[th] century, it can only be traced through isolated remarks. In the second half of the 19[th] century, authorship disputes, in particular the Bacon-Shakespeare controversy, fueled more systematic research. A detailed history of quantitative authorship attribution appears in (Grieve 2005). We shall complete it, presenting earlier sources on the proportion of vowels to consonants.

On the history of Statistics, we shall refer to (Hacking 1990), and (Stigler 1999). A major change took place in the perception of statistical concepts, during the second half of the 19[th] century: partly inspired by Darwin's theory of evolution and the controversy around Quetelet's "Average Man", statistics gradually moved from the mechanist perception inherited from astronomy, of a ʻtrue valueʼ measured with random errors, to the modern notion of intrinsic variability. Letter frequency distributions were regarded by early statisticians, in particular Babbage and Quetelet, as one of the "Constants of Nature and Arts", just as they had always been viewed by cryptologists and linguists. The perception gradually evolved to the modern notion of an information source defined as a stochastic process, through the works of Markov, then Kullback and Shannon: Cryptanalysis was crucial in the creation of information theory. At the same period, Zipf's observations on the relative frequency of words initiated a seemingly different corpus of studies in Statistics. Both converged in the second half of the 20[th] century through the works of Mandelbrot (Brillouin 2004).

A short history of letter counting will first be presented. Then the three domains, Cryptology, Quantitative Linguistics, and Statistics, will be examined in that order, focusing on the interactions with the other two fields through letter counting. As a conclusion, the eclectism of past centuries scholars, a generalized background in humanities, and the familiarity with cryptograms, will be identified as contributing factors to the mutual enrichment process which is described here.

**Letter counting: from rope-dancing to industrial mechanics**

When Yule referred to text counting as being "of the most venerable antiquity and august associations", he meant the Jewish tradition of the Masorah, about which he cited Hasting's "Dictionary of the Bible" (Yule 1944, 7, see also Kelley et al. 1998)*:*

> They *counted* the verses and the words of each of the 24 books and of many sections; they reckoned which was the middle verse and the middle word of each book; nay, they counted the letters both of particular sections and even of whole books. Thus they could specify the middle word and the middle letter in the Torah, the middle verse and the middle letter of the Psalms. They counted also the frequency of occurrence of words, phrases, or forms, both in the whole Bible and in parts of it.

Thus, for proofreading or salary purposes, Bible copyists have counted letters and words for many centuries. Yet, even before the Masoretes, long before any text had ever been counted, evidence was there that letter frequencies could be biased at will. Some authors did not fail to see it. Here is what ʿAlī ibn ad-Durayhim (1312–1361) remarked (Mrayati et al. 2003b, 100).

> This is the letter order of frequency in the Holy Koran, although this order may differ in other language usages. Some deliberately encipher poetry and prose, dispensing with the letter [ā], or without letter-dotting or without idle particles [those that do not affect the parsing of what follows].

Dispensing with one or several letters in writing is the principle of *lipograms*, already practiced by some Greek poets (Perec 1973). Georges Perec, being himself the author of an impressive novel of about 360,000 letters that does not contain a single *e*, had perfectly seen the relation between lipograms, cryptology, and the Holy Scriptures: "the Book is a cryptogram of which the Alphabet is the cipher" (Perec 1973, 73). Since by lipograms or any other form of constrained writing, letter frequencies can be skewed, why then did all cryptologists, linguists, and statisticians fail to mention it before the end of the 19$^{th}$ century? Part of the explanation may be that word games, even practiced by the best authors, were not taken very seriously before the 20$^{th}$ century. In his book on anagrams, Henry Benjamin Wheatley devotes 3 pages to the history of lipograms (Wheatley 1862, 21–23). Here is how they begin.

> Lipograms are perhaps the most ridiculous of those follies; it is of them that De Quincey speaks when he says "Recalling that rope-dancing feat of some verse-writers, who through each several stanza in its turn, had gloried in dispensing with some one separate consonant, some vowel or some diphtong, and thus achieving a triumph such as crowns with laurel that pedestrian athlete who wins a race by hopping on one leg, or wins it under the inhuman condition of confining both legs within a sack."

After the art of printing was invented, letter counting continued, for technical reasons: typographers needed to know what quantities of each character had to be prepared. In 1828, Adolphe Quetelet (1796–1874) is not yet the founding father of statistics in social sciences. When he publishes in his new journal, an article entitled "On the ratios between the number of letters that compose the alphabet in different languages", he does so under the section "Industrial Mechanics", and not Statistics as one would have expected from him (Quetelet 1828, 339–341). Among his motivations, typography clearly comes before cryptology:

> It is important for typographers and foundry workers, to know the proportions of letters that compose the alphabet in the different languages, in order to form their stock. The knowledge of these proportions may also interest the philologist; on the other hand, suspicious governments have sometimes used it to decipher correspondences in which the letters of the alphabet were substituted one for the other. Since these proportions are generally little known, we have thought appropriate to record here those of three alphabets [Dutch, French, Italian].

In the following volume, Quetelet renders a detailed account of a thesis defended by Marcel Hayez, from the family of typographers who publishes the "Correspondance Mathématique et Physique" (Quetelet 1829, 391–394). Hayez has provided detailed letter and character counts for French, Latin, Italian, Portuguese, Greek, English, Dutch, and German, over corpuses of several thousand letters; quite naturally, typography is both his motivation and his application domain.

That typographers might have a precise idea of character frequencies, came to the mind of Samuel Morse in 1838 when he was in search for a binary code that would optimize time, hence match the simplest codes to the most frequent letters. Instead of painstakingly counting letters by himself, he simply walked to a printing office and noted the quantities that he found there (Morse 1914, 69). Both letter frequencies and the Morse alphabet later contributed shaping the QWERTY keyboard, that we are still using (Yasuoka 2011). Contrarily to Morse, it seems that Alfred Butts, the inventor of the Scrabble, counted letter frequencies from newspapers front pages, to decide the value and number of the letters in his new game, at least if web resources are to be believed on that point.

**Cryptology: don't spare either labor or paper**

According to usage, we distinguish *cryptography* (ciphering) from its contrary *cryptanalysis* (deciphering). Not all Cryptology treatises deal with cryptanalysis. Those who do, usually begin with a section on letter frequencies. Some quotations, given in chronological order, will help understand how the frequency analysis method has been perceived by cryptographers along the ages. A very precise and rational procedure for the Arabs, it came to be one of the tools of a largely divinatory art among westerners of the 15$^{th}$ to 18$^{th}$ centuries.

The clearest exposition of the method may be the first one, by Al Kindī (Mrayati et al. 2002, 124):

> Among the expedients we use in cryptanalysing a cryptogram if the language is already known, is to acquire a fairly long plaintext in that language, and count the number of each of its letters. We mark the most frequent letter ʿfirstʾ, the second most frequent ʿsecondʾ, and the following one ʿthirdʾ, and so forth until we have covered all its letters. Then we go back to the message we want to cryptanalyse, and classify the different symbols, searching for the most frequent symbol of the cryptogram and we regard it as being the same letter we have marked ʿfirstʾ in the plaintext; then we go to the second frequent letter and consider it as being the same letter we have termed ʿsecondʾ, and the following one ʿthirdʾ, and so on until we exhaust all the symbols used in this cryptogram sought for cryptanalysis.

Indeed, Al Kindī had based his counts on 7 pages totaling 3667 letters, a feat that westerners were not to accomplish for a very long time. Three centuries later, ʿAfīf ibn ʿAdlān (1187–1268) proposed a classification of letters into three frequency classes, and introduced frequencies of digrams (Mrayati et al 2003a). His contemporary ʿIbrāhīm ibn Dunaynīr (1187–1229) repeated Al Kindī's counts, and could only confirm: "So it came home to me the validity of the statement of Yaʿkūb Al Kindī, peace be on his soul." Here is how he explains the method (Mrayati et al. 2005,

84).

> I say: The type of encipherment characterized by changing the forms of letters is achieved by devising shapes or symbols not attributed to letters at all. In this method every letter is represented by a symbol that is unique to it. The cryptanalysis is accomplished by counting the symbols of the cipher message, and establishing the frequency of occurrence for each symbol, by affixing the frequency number to the respective symbols of the cipher. Having done that, you dispose the symbols in order of frequency precedence, designating the most frequently-occurring symbol in its locations throughout the cryptogram. Do the same with the next frequently-occurring symbol, and so forth until you exhaust all the symbols of enciphered letters. Now place the highest-frequency symbol against the highest frequency letter of the Arabic alphabet, doing the same with the rest, conformably with their order of frequency. Keep going in the same vein until you use up all the letters and symbols.

Some Arab treatises went further than Al Kindī on particular linguistic or statistical points; but from the point of view of Cryptology, the essential was there. Confronted to the same problem, cryptographers of the Italian Renaissance expressed their findings in a different way. The first of them is Alberti, already cited in the introduction. Apart from a very interesting linguistic remark that will be examined in the following section, Alberti does not give exact frequencies, but rather indications of relative orderings, and mutual relations between letters. After a careful (but not quantified) linguistic examination of Latin, he somewhat maliciously concludes (Alberti 2010, 177):

> We have discussed the number and the order of the letters and those vowels that are most frequently encountered in writing, of the possible association of vowels with regards to the consonants, we have said which are most frequently encountered and regarding their arrangement, how they behave in addition or in sequence. From what has been said thus far you will quite easily deduce (if I am not mistaken), how astute minds will perceive a way to decipher and interpret coded writing in a most easy way.

A few years after Alberti, when Cicco Simonetta (1410–1480) exposes his method for deciphering, his first seven rules aim at deciding whether the enciphered text is in Latin or "Volgare." In the other 6 rules, he does give observations to isolate vowels, and recognize the other letters, but never mentions frequencies (Perret 1890). Italian cities, and in particular the Vatican, are known to have had excellent cryptanalysts for a long time; among them the Argenti family. After his dismissal in 1605, Matteo Argenti compiles an impressive manual, exposing his method of cryptanalysis quite clearly. Sixty two rules are listed; only the 22[nd] mentions that frequent letters are likely to be vowels (Meister 1906, 166).

At about the same period, François Viète (1540–1603) served as the decipherer of the French court. His notes, written shortly before his death, prove that he used frequency analysis, even though it was only one among several other tools (Pesic 1992, 13):

> One must note all the sorts of figures, whether ciphers or jargon, and count how many times they occur, then note all the sorts of figures which precede or which follow and compare the most frequent in order to discover the same words, and the same meanings. Don't spare either labor or paper.

By the end of the 17th century, the frequency analysis method was well understood by everyone in the trade, and even popularized in treatises, such as John Falconer's (Falconer 1685, 8–10).

> To proceed regularly therein, you must endeavour to *distinguish betwixt the Vowels and Consonants*. And first, the Vowels generally discover themselves by their frequency; for, because they are but few in number, and no word made up without some of them, they must frequently be used in any Writing : however, it may by accident fall out, that some of the Consonants shall be oftener found in an Epistle.
> […]
> *To distinguish one Vowel from another*, after you have made the most probable Suppositions in separating them from the Consonants, 1. Compare their frequency, and *e*, as we observed before, is generally of most use in the *English* Tongue, next *o*, then *a* and *i*; but u and *y* are not so frequently used as some of the Consonants.
> […]
> *To distinguish one Consonant from another*, you must 1. (as before) observe their frequency. Those of most use in *English* are *d*, *h*, *n*, *r*, *s*, *t*, and next to those may be reckoned *c*, *f*, *g*, *l*, *m*, *w*; in a third rank may be placed *b*, *k*, *p*; and lastly, *q*, *x*, *z*.

David Arnold Conradus, published in 1732 a Latin book *Cryptographia denudata, sive Ars Decifrandi,* (Conradus 1732), of which he gave, ten years later, a more widely accessible version in a series of articles published by "The Gentleman's Magazine" (Conradus 1742, 134).

> *Prop 2*] Every Language has, besides the Form of its Characters, something peculiar in the Place, Order, Combination, Frequency, and Number of the Letters.
> *Rule 1*] In Deciphering, Regard is to be had to the place, Order, Combination, Frequency, and number of the Letters
> 　　　[…]
> *Prop. 6*] The Vowels, generally five, are four times out-numbered by the Consonants, the Vowels must therefore recur most frequently.
> *Rule 6*] The Letters that recur most often are Vowels.

To the best of our knowledge, only one *quantified* assertion about frequencies can be found in a Cryptology treatise between Alberti and the 19th century: it appears in Christian Breithaupt's "Ars Decifratoria" (Breithaupt 1737, 109); it is much less precise than Alberti's observations.

> The frequency of letters should be noted in general, since in any language vowels are more numerous than consonants. The reason for making these observations is obvious. Actually, for a given number of vowels, the corresponding number of consonants must be larger by five fourths; it cannot be otherwise, vowels being more frequent than consonants.

Half a century after Conradus and Breithaupt, Philip Thicknesse (1719–1792) describes again the same deciphering method, with hardly any novelty (Thicknesse 1772, 22–25).

> The next thing to observe is, what letters occur oftenest, and those you may conclude are vowels, and that which is most frequent, to be an *e* –as *e* in English occurs oftener than any other letter.
> 　　　[…]
> Observe also, that *i*, in English, never terminates a word, nor *a* or *u* except in *sea*,

> *you*; or those; and thus by comparing the frequency of the letters, you will generally find *e* occur the oftenest: next *o*, then *a*, and *i*; but *u* and *y*, are not so often used as some of the consonants. Among the vowels, *e* and *o* are often doubled; the rest scarce ever: and *e* and *y*, often terminate words, but *y* is less frequent, and consequently easily distinguished.
> To find out one consonant from another, you must also observe their frequency, *d*, *h*, *n*, *r*, *s*, *t*, and next to those, *c*, *f*, *g*, *l*, *m*, *w*, in a third rank may be placed *b*, *k*, *p*, and lastly *q*, *x*, *z*. This remark, however, belongs to English; for in Latin the consonants are *l*, *r*, *s*, *t*; next *c*, *f*, *m*, *n*; then *d*, *g*, *h*, *p*, *q*; and lastly, *b*, *x*, z.

Admittedly, before computers, counting thousands of letters to establish a frequency table was a time consuming and not very interesting activity. Most 20[th] century manuals contain such tables, but rare are those who do not reproduce earlier counts. For instance, in the French edition (Sacco 1951) of one of the most translated and re-published cryptology manuals, the source of the tables is cited as (Valerio 1893). We have seen (Klüber 1809) cited as a source of tables, yet no tables were available in the edition that we could consult. The first European manual of Cryptology to give frequency tables for different languages, that we are aware of, is Charles François Vesin's "Traité d'Obscurigraphie" (Vesin 1838, 100–104). As we have seen, similar tables had already been established earlier, by Quetelet and Hayez. Vesin only announces "having tried a corpus of one thousand letters." The table comes with some warning: "There may doubtlessly be some difference in the distribution of letters, between some period of writing and another period; but this difference will never be exorbitant." By the turn of the 20[th] century, such warnings became more systematic. Yet, their authors still firmly believed in the stability of letter frequencies for ʿnormalʾ writing. Author of a "Manual for the solution of military cipher", intended for Army Service Schools shortly before the United States entered WWI, Parker Hitt considers that deciphering offices "should be provided with tables of frequency of the language of the enemy, covering letters and digraphs." Explaining what he means by "tables of frequency", he appears aware of the stylistic problem, which he downplays (Hitt 1916, 4).

> In fact, if ten thousand consecutive letters of a text be counted and the frequency of occurrence of each letter be noted, the numbers found will be practically identical with hose obtained from another text of ten thousand letters in the same language. The relative proportion of occurrence of the various letters will also hold approximately for even very short texts.
>     Such a count of a large number of letters, when it is put in the form of a table, is known as a frequency table. Every language has its own frequency table and, for any language, the frequency table is almost as fixed as the alphabet. There are minor differences in frequency tables prepared from texts on special subjects. For example, if the text be newspaper matter, the frequency table will differ slightly from one prepared from military orders and will also differ slightly from one prepared from telegraph messages. But these differences are very slight as compared with the differences between tables of two different languages.

Hitt's presentation is still quite far from the modern view of frequency tables. Here is how they are described in a present day manual (Bauer 2006, Chapter 15).

> Moreover, even long texts normally show considerable fluctuations of character frequencies. ʿTheʾ frequency distribution of English is a fiction, and at best the military, diplomatic, commercial, or literary sub-languages show some homogeneity; indeed even the same person may speak a different language on the circumstances. Correspondingly, statistics on letter frequencies in different languages are quite variable. Moreover, most of the older counts were based on texts of only 10,000 or fewer letters. For the frequency ordering there are already great differences in the literature.

After having traced the frequency analysis method through cryptology manuals, several questions remain to be answered. The first one is: why was the method considered useful for such a long time? Many manuals, after having exposed the vulnerability of monoalphabetic substitution, propose alternatives. Alberti had already described a polyalphabetic cipher, soon improved by Tithemius, Bellaso, and Vigenère (Kahn 1996, Mollin 2005). Why weren't these ciphers systematically used from them on? In 1839–1841, Edgar A. Poe earned himself a high reputation as a cryptanalyst, by solving simple substitution ciphers. However, he does not seem to have had much more competence than he had found in (Conradus 1742): see (Friedman 1936), and (Poe & Peithman 1986, 283). According to (Kahn 1996), codes that were vulnerable to frequency analysis remained in use throughout WWII. The explanation usually retained is that polyalphabetic ciphers were too complicated to use, in particular by military services. Another explanation may be diplomatic secrecy. From the 15$^{th}$ century on, most courts used their own ciphers, usually more complicated than monoalphabetic substitution, yet still vulnerable to frequency analysis; but no one wanted the others to know that they could solve their codes. This will be illustrated through a well known anecdote about Viète's feats as a cryptanalyst in the diplomatic war between the kings of France and Spain, Henry (Henri IV) and Philip (Felipe II). Here is how it is reported in (Kahn 1996, 82).

> Meanwhile, Philip had learned, from his own interceptions of French letters, that Viète had broken a cipher that the Spanish –who apparently knew little about cryptanalysis– had thought unbreakable. It irritated him, and thinking that he would cause trouble for the French at no cost to himself, told the pope that Henry could have read his ciphers only by black magic. But the tactic boomeranged. The pope, cognizant of the ability of his own cryptologist, Giovanni Batista Argenti, and perhaps aware that papal cryptanalysts had themselves solved one of Philip's ciphers 30 years before, did nothing about the Spaniards' complaint; all Philip got for his effort was the ridicule and derision of everyone who heard about it.

The historical context is reported by (Pesic 1997) in similar terms: "Spain seems not no have had cryptanalysts, which might explain Philip's attitude." Complacently spread by many French sources along the centuries, that story resembles more propaganda than historical truth. At the same period, Philip also had his own cryptanalyst, Luis Valle de la Cerda, who was as much able to decipher Henri's dispatches as Viète was with Philip's (Carnicer Garcia & Rivas 2005, 73). Contrarily to Henry, Philip probably preferred the French to have some amusement at his own expense, rather than letting them know about Valle de la Cerda's achievements.

Secrecy may not have been the only reason for the relative lack of popularity of frequency analysis among cryptologists. It must be emphasized that the method was neither necessary nor sufficient to solve the ciphers cryptanalysts were usually confronted with. That it was not necessary has been illustrated in particular by Simonetta and Argenti: there are many other ways to guess vowels, by looking for isolated characters, word endings, digrams, and trigrams. The method was not sufficient for two reasons. One is the length of encrypted messages: many authors have pointed out that the shorter the message, the less reliable the order of frequencies. That point will be examined in the Statistics section. The other reason is that, even when a monoalphabetic code had been chosen, plenty of ways to level out frequency differences were available, and recommended in manuals; the simplest was to alternate different characters for the same letter. Using or not frequency analysis, actual deciphering was an arduous task, requiring at least as much intuition as method; hence Viète's exhortation: "don't spare either labor or paper!"

**Quantitative Linguistics: if sholars knew the law of averages as well as mathematicians…**

Reading the collection of the "Arabic Origins of Cryptology" (Mrayati el al. 2002 ff.), one cannot but remark how Arab authors often went much further into the analysis of the language, than strictly needed for cryptanalysis technical purposes: their treatises contain many profound remarks on the structure of language, on the alternation of vowels and consonants, on the differences between poetry and prose, on rhythm and prosody. The symbiosis between Arab Cryptology and Language Sciences is described in (Mrayati et al. 2002, 46–53). As stated there, the relationship "was not unilateral. Biographical sources confirm that the leading scholars of language were also involved in cryptology." The same can be said of many western scholars. Our purpose here is not to analyze the relationship between Cryptology and Linguistics in general, but rather to describe, through some particular examples, how some linguistic analyses of the language, involving letter counts, may have fostered the emergence of Quantitative Linguistics as a new field. Historical accounts of Quantitative Linguistics written so far have given the main role to autorship attribution problems (e.g. Williams 1956, Bailey 1969, Peng & Hengartner 2002). This can be explained on the one hand by the hotly debated Bacon-Shakespeare question (Stopes 1888), on the other hand by the influential work of (Yule 1944). We believe that the history of the field should be re-evaluated, taking into account earlier discussions on the vowel-consonant proportion, three examples of which are given in this section.

A remark in ad-Durayhim's treatise (Mrayati et al. 2003b, 100), already cited, shows that he was aware of possible differences in letter frequencies, under a change of style. Yet we have not seen any such remark being *quantified* in Arab treatises. To the best of our knowledge, the first stylometric observation in history is due to Alberti, in "De Componendis Cifris" (Alberti 2010, 173, Meister 1906, 127).

> From my calculations, it turns out that in the case of poetry, the number of consonants exceeds the number of vowels by no more than an octave, while in the case of prose the consonants do not usually exceed the vowels by a ratio greater than a sesquialtera. If in fact we add up all the vowels on a page, let's say there are three hundred, the overall sum of the consonants will be four hundred.

That remark has been analyzed in (Ycart 2012). The interpretation proposed there is the following. The proportion of vowels for poets may be denoted by P, being aware that what precisely is meant by ʽpoetsʼ and ʽvowelsʼ, is a matter of convention. Alberti's assertion can be mathematically translated into $(1–P)–P<1/8$, or else $P>7/16$. Alberti opposes poets to "rhetores", i.e. orators, for which the proportion of vowels against consonants is said to be "sesquitertia" in the Latin text, or else four against three. If R denotes the proportion of vowels for orators (again a matter of convention), what is stated is $(1–R)/R < 4/3$, or else $R>3/7$. The difference between $7/16=43.75\%$ and $3/7=42.86\%$ is smaller than 1%. Yet, using a corpus of Latin texts, from poets and orators, it could be proved that there was indeed a significant difference in the use of vowels, and that Alberti could have detected that difference by letter counts. Moreover, the last part of the citation gives an idea on his sample sizes: mentioning a page of 300 vowels and 400 consonants might have been a way of indicating that his observations were supported by counts on large enough sets of letters. Because of that remark, some have seen in Alberti a pioneer of ergodic and information theories (Mercanti & Landra 2007, note 14). This may be exaggerated; notwithstanding, the two sentences quoted above are quite singular. Alberti himself, while commenting at length vowels, consonants, their relative frequencies and positions, never gives any other number. As already stated, the only other quantified remark found in a Cryptology treatise, that of Breithaupt, falls short of Alberti's observation as far as accuracy is concerned.

Our next example does not come from Cryptology, yet it deals with the same problem of vowel-consonant proportion. William Owen Pughe (1759–1835) a grammarian and lexicographer, devoted his life to the illustration and defense of the Welsh language, producing in particular a grammar, and a monumental Welsh-English dictionary. Between 1796 and 1816 at least, he regularly contributed to "The Monthly Magazine", under the pen name of "Meirion", sending letters on etymology, poetry, "the affinity of Welsh and Hebrew", and so on (Davies 2002, 36). On January 1st, 1799, he wrote the following (Pughe 1799).

> In looking into the new edition of Chamber's Cyclopaedia, sometime ago, I casually met with a remark upon a subject, which had relation to language, wherein the *Welsh* and the *Dutch* were pointed out, as abounding more with consonants, than most, if not all the European tongues. I well knew that such a statement was proverbial, as a vulgar prejudice: but I became a little angry, at finding it had obtained a place in one of the first philosophical dictionaries of the present age; and, not being able to efface the impression from my mind, I had recourse to the finding a tolerably exact arithmetical certainty, as to the fallacy of such an observation. The method, adopted as the most eligible, was to fix upon the mean number of vowels to a hundred consonants, in different languages; and to exhibit the results in a table. As the conclusion, to be drawn from it, tends to establish a point, if not of importance, at least of some curiosity, you may be induced, Sir, to give it insertion, in your valuable repository.
> [table for Welsh, Greek, Spanish, Italian, Latin, French, German, Dutch, English]
> As the French, and the English, differ so considerably, in pronunciation, from what they appear in orthography, the following comparison shews the reduced numbers of the vowels and the consonants[…] The English is very variable, with respect to the proportion of vowels and consonants: that of the consonants is much greater in the scripture style, than in elegant writing, and more especially that which is scientific, from its containing more words derived from the learned languages. In the bible, the compass of the variation, in the number of vowels, is generally from about 68 to 50; but the medium may be settled at 56 to a 100 consonants. In polished writing, the medium number of vowels may be settled at 66; and the mean between the two styles will be 61, the number inserted in the foregoing table.
>
> The compass of variation in the Greek is considerable. I have found 150 vowels to 100 consonants; and frequently as low as 86. The other languages are pretty close to the average number, given in the table; the Welsh seldom deviates three vowels from the mean number.

By modern standards, some scientific objectivity might be wanting. Yet, with Statistics still in its infancy, Pughes shows a quite remarkable care in assessing the variability of his counts and find a "tolerably exact arithmetical certainty": recall that by 1799, the word "Statistics" had just been given its modern meaning by Sir John Sinclair a few years before (Hacking 1990, 26). Counting the number of vowels per one hundred consonants is unusual, but perfectly correct if enough samples are taken to assess the "medium" number of vowels, together with a reasonable "compass of variation." The observation that the ratio vowel-consonant depends not only on the language, but also on the style is there, and also the distinction between pronunciation and orthography. In his defense of Welsh, Pughe means to prove that it has more vowels than the other tongues, and he accepts (and even supports) the idea that "the harmoniousness of a language depends much upon the proportion of the vowels to the consonants." The sentence comes from the article "Vowel" of the Encyclopedia Americana, and expresses what was a commonplace prejudice in the beginnings of comparative linguistics. It was written by Francis Lieber (1800–1872), an early specialist of linguistic relativism (Andresen 1996, 114–119). We shall quote from his unusually long article "Consonants" of the same Encyclopedia Americana (Lieber 1835). There, he exposes a statistical study of his, comparing the ratio vowel-consonant among classical poets of different languages.

> It occurred to the writer, while preparing this article, that it would lead to interesting results, if the proportion of the vowels and consonants, in the different languages, could be ascertained; but the conclusions, to which he has been led by such investigation as he has bestowed on the subject, are rather to be regarded as indications of what might be learned from more thorough inquiries, than as facts from which general deductions can be safely drawn. In making the comparison, passages have been taken from the popular poets of different countries. The different passages were in the same measure, or in measures very similar, so that the number of syllables in each would be very nearly the same.

Of course, Lieber distinguishes the "orthographic proportion" from the "phonic proportion." He also knows that euphony cannot only come from the frequency of vowels. Interestingly, he is quite cautious on methodology issues and aware of the necessity of a rigorously selected corpus of sufficient size.

> To give any thing like accuracy to such investigations, it is obvious that the results ought to be taken both from prose and poetry, also from many different writers, and the language of conversation.

In the conclusion of the article, Lieber clearly expresses his ambition: nothing less than a new branch of Linguistics, that would quantify the elements composing a language, exactly as the then newly founded stoechiometry quantifies atoms and molecules in Chemistry: quite an interesting comparison for a German born political philosopher, trained in Mathematics at the University of Jena (see Hufbauer 1982).

> It is easily seen, that, in the languages of Latin origin, the proportion of consonants to vowels is much smaller than in the Teutonic idioms. To compare the proportions of consonants to vowels, in such different families of languages; to show the proportions of the gutturals, labials, &c., of the different idioms; and again, the proportion of these letters in the various families of languages, or according to the different parts of the earth to which they belong, as Asiatic, European, &c. languages, and many other calculations – might lead to very interesting conclusions. This branch of philology might be compared to the new department of *stoechiometry* in chemistry, which treats the proportions of the quantities of the elements in a state of neutralization or solution – a branch of science which every day becomes more important, and which has been illustrated by the labors, past and present, of a Berzelius, Klaproth, Döbereiner and others.

Lieber's manifesto came 16 years before another declaration of intent, much less precise, yet most often cited as a birth date for Statistical Linguistics (Bailey 1969, Grieve 2005): a letter to one of his friends by the famous mathematician Augustus De Morgan (1806–1871), dated August 18, 1851, where he suggested statistical counting to settle authorship disputes (De Morgan 1882). De Morgan did not mean to go through the tedious process of counting by himself:

> It has always run in my head that a little expenditure of money would settle questions of authorship in this way. The best mode of explaining what I would try will be to put down the results I would expect as if I had tried them.

In the conclusion of his letter, De Morgan expresses his absolute trust in the law of large numbers: he is perfectly acquainted with the statistical developments of his time, but the paradigmatic shift of statistics had not yet taken place (Hacking 1990).

> If scholars knew the law of averages as well as mathematicians, it would be easy to raise a few hundred pounds to try this experiment on a grand scale. I would have Greek, Latin, and English tried, and I should expect to find that one man writing on two different subjects agrees more nearly with himself than two different men writing on the same subject. Some of these days, spurious writings will be detected by this test.

A complete account of how De Morgan's ideas, published after his death, were implemented by Thomas C. Mendenhall (1841–1923), and eventually developed into a full methodology for quantitative authorship attribution, has been given by (Grieves 2005). It will not be reproduced here. The interested reader is also referred to (Yule 1944), already cited. It is quite instructive to compare the idealistic confidence of De Morgan, who did not actually count letters, and wrote before the great statistical shift, to the wariness of Yule, who was one of the very best statisticians of his time, and had counted himself many texts (Yule 1944, preface).

> When I had advanced some way in that particular study, it became only too clear into thorny a field I had strayed. Statistics of literary vocabulary proved to have their own special problems, their own peculiar difficulties and sources of fallacy, which no one apparently had made any attempt systematically to explore.

We conclude this section with a quite different contribution of Quantitative Linguistics to Statistics. Back in the 10$^{th}$ century, ibn Wahab al-Kātib had the following interesting remark on Arab phonology (Mrayati et al 2007, 110):

> The closer the outlets of two letters, the heavier on the tongue in articulation. The Arabs tend by nature to use what is inherently light and easy to pronounce, and to avoid using all that is awkward. It follows that they hardly take up two letters of the same outlet or of two equal outlets. If, however, such letters happen to occur together, they tend to assimilate the one into the other.

That remark was not to be quantified until ten centuries later, with Zipf's work on the "law of least effort" (Zipf 1949). Nowadays, statistical applications of Zipf laws and the like (Benford law, Pareto law) are ubiquitous, from econometrics to physics, biology, and computer science (Reed 2001).

**Statistics will reach an accuracy level…**

It cannot be asserted that letter counting has always been a central preoccupation to statisticians. Yet it was present at several crucial stages in the development of the field: the Arab beginnings, the organization into a new scientific discipline in the middle of the 19$^{th}$ century, the early development of stochastic processes by Markov, and the creation of information theory by the cryptanalysts of WWII. These four periods will be examined successively.

As we have seen before, Arab cryptologists have made quite substantial letter countings (3667 letters for Al Kindī). They certainly deserve to be regarded as pioneers of Statistics, at least in the sense of data collection; presenting them as specialists of data *treatment* would be an overstatement. However, it must be noted that they did have a clear idea of variability and sample size. This will be illustrated by two quotations from ʿAfīf ibn ʿAdlān's treatise. The first one shows how he expressed letter frequencies, comparing them to the most frequent. Notice the last sentence of the passage (Mrayati et al. 2003a, 48).

> As far as frequency of occurrence is concerned, some Arabic letters are abundant, and these are seven in number, i.e.: [ā], [l], [m], [h], [w], [y], and [n]. If the letter [ā] occurred in a piece of writing as many times as, say, 600 times, it has been found as a corollary, that the letter [l] would occur about 400 times, add or take few, the letter [m] 320 times, the letter [h] 270 times, [w] 260 times, [y] 250 times, and the letter [n] 220 times. This is so in all likelihood, but there may be some variation in this order.

The second quotation from the same treatise gives a quantified indication on the length of the text (Mrayati et al. 2003a, 52).

> The length of the text to be cryptanalysed should be at least in the neighborhood of ninety letters as a rule of thumb, because the letters thus would have had three rotations. Yet, the number of letters may be less than that in certain cases.

Arguably, the organization of Statistics into a new branch of Science can be dated to a visit of Adolphe Quetelet to London in 1933, and his meeting with several statisticians there (Hacking 1990, 61, Babbage 1864, 434, Quetelet 1873, 156). Among them, Charles Babbage (1781–1871). Quetelet and him had made friends at a dinner in Paris, in 1826 (Quetelet 1873, 152), and they had regularly corresponded since. An interesting character with multiple interests, Babbage is best known nowadays for his "Analytical Engine", the ancestor of computers. He also deserves a place in the histories of Cryptology and Statistics. In Cryptology, he is credited for being the first to have cracked the famous Vigenère code (Mollin 1990, 74). His own feelings were ambivalent: "Deciphering is, in my opinion, one of the most fascinating of arts, and I fear I have wasted upon it more time than it deserves" (Babbage 1864, 235). During that "wasted time", he acquired a solid experience of lexicography, counting and sorting letters and words. In Statistics, Babbage was a passionate supporter of the "Law of Averages." He proposed to gather in a huge set of data the "Constants of Nature and Arts", that he had grouped in 19 chapters. Among these constants, a good number were of statistical nature (Hacking 1990, 52). Chapter 18 of the list certainly owes much to his experience as a cryptanalyst (Babbage 1832, 338).

> 18. Tables of the frequency of occurrence of the various letters of the alphabet in different languages, – of the frequency of occurrences of the various letters at the beginnings or endings of words, – as the second or as the penultimate letters of words, – of the number of double letters occurring in different languages, – of the proportion of letters commencing surnames amongst different nations.

It was more than a vague intention: on the same year Babbage sent to Quetelet a frequency table of double letters in 10,000 words of 7 different languages (Quetelet 1832, 135–136). Twenty years later, during the preparation of the first international congress, he insists again (Quetelet 1873, 162):

> On reading the printed views of your statistical congress, it appears to me that the

> constants of nature are within your limits, and I am more than ever convinced that even an attempt to collect a portion of them would contribute largely to the advancement of knowledge and to the economy of the time of its cultivators.

Babbage devoted one of his talks of 1853 in Brussels to the Constants of Nature and Arts, and Quetelet agreed to reproduce (Babbage 1832) in the proceedings (Quetelet 1853, 226). Quetelet fully supported Babbage's views on the law of averages: his own theory of the "Average Man" had been conceived before 1831 (Stigler 1999, 51–65). As already said, Quetelet had himself counted letters in 1828. Yet, he never was really convinced. Forty years later, he remembered the episode more like one of Babbage's quirks, than a true contribution to the Advancement of Knowledge.

> We have already mentioned, at the beginning of this article, the singular tendency that our friend had, to deal more specially with questions that others lost sight of. I thought appropriate not to conceal from him a remark that I had already made in this respect; presumably, my observations made him smile, and he transmitted to me a new writing on the subject. That time, he had lost sight of high speculations of political economy: he had dealt with researches that he had made on the number of times that letters are doubled in a word, for instance on ten thousand different words. It may be curious to see the results that he got, comparing the English, Italian, German, and Latin languages. But what advantage can be expected in general from such a study?

Admittedly, direct profit for the happiness of mankind was remote; but one cannot help wondering whether by counting more letters, Quetelet and Babbage would have realized earlier, how variable some "Constants of Nature and Arts" can be.

Neither of them seems to have been aware of the work of Viktor Yakovlevitch Bunyakovsky (1804–1889). The contributions of this Russian mathematician to Statistics deserve a more important place than they have been given so far, in particular in connection to our subject. In 1847, Bunyakovsky published in a literary journal an article entitled "On the possibility of introducing definite measures of confidence in the results of some observation sciences, mainly statistics" (Petruszewycz 1979b, 47). As the title says, confidence intervals are rigorously defined, and Bunyakovsky concludes with a prophetic vision of the discipline.

> One can nearly predict that, some time from now and perhaps very soon, Statistics, following observation sciences, will reach an accuracy level that, while remaining by essence below that of astronomic results, will still be way above what it is in the present state of that science.

As Bunyakovsky himself says, he could have terminated his article there. But he goes on to present a new application.

> The new application relates to grammatical and etymological research, as well as comparative philology. At first glance, such investigations seem alien to mathematical analysis, nevertheless it can be said with certainty that a huge field for rigorous mathematical applications lies there. My assertion is not based upon more or less rigorous assumptions and conjectures, but upon a critical examination of the discipline, upon some trials that I have already realized, and upon analytical formulas that I have deduced to define numerical probabilities of various linguistic derivations.

Bunyakovsky then goes on, giving examples of assertions that could be submitted to rigorous statistical proof, detailing data that should be collected on languages, and finally, advocating a collaboration between mathematicians and philologists. He also mentions the possibility to publish his own research on the subject, "should another opportunity arise"; unfortunately, no such work of his has been found to this date. Nevertheless, he not only paved the way for rigorous assertions in statistics, grounded on precise confidence intervals, but he also set up an even more ambitious program than that of Lieber: Bunyakovsky deserves to be counted as well among the founders of Quantitative Linguistics.

When Andrei Andreievich Markov (1856–1922) began working on his chains of events in 1906, his motivation was mathematical: he wanted to prove that stochastic independence is not a necessary condition for the law of large numbers. Yet the only illustration that he ever gave was again the proportion of vowels to consonants: he counted 20,000 letters in Pushkin's "Eugen Onegin", 100,000 letters in Aksakov's "Childhood years of Bagrov's grandson", and made a very thorough statistical study of his counts, proving in particular that vowels and consonants did not alternate independently. The mathematical motivations of Markov, the connection of his work with cryptology and linguistic studies of his time, have all been amply documented, by Micheline Petruszewycz (Petruszewycz 1979a, 1983), and more recently by David Link (Link 2006a, 2006b) who also provided an English translation of Markov's famous article on Eugen Onegin, together with Gloria Custance (Markov 2006). The reason why Markov chose his linguistic application has also been thoroughly discussed: Micheline Petruszewycz convincingly argues that Bunyakovski's article may have been a source of inspiration (Petruszewycz 1979b). The obvious vowel-consonant filiation from Alberti to Markov through Pughe and Lieber, should not hide the major shift inaugurated by Markov in the mathematical modeling of language. From Markov on, a succession of letters would be seen as a dynamic source of information, i.e. as a *stochastic process*. The approach is perfectly explained in (Shannon 1948), the founding article of information theory by Claude E. Shannon (1916–2001). There, Shannon proposes several models of increasing likelihood and information content, for the random generation of texts. He explicitly refers to Markov chains, and cites Maurice Fréchet (1878–1973), who helped popularize Markov's ideas in the west (Link 2006b). Shannon later detailed the application of his information theory to Cryptology (Shannon 1949) and Linguistics (Shannon 1951). The history of information theory has been accounted for in several works: (e.g. Pierce 1973, Brillouin 2004, Cover & Thomas 2006). It would be simplistic to reduce this major landmark of Statistics to the single influence of Cryptology. However, the deep involvement of some of the major actors into the war effort during WWII cannot be downplayed (Kahn 1996, Mollin 2005). Not only Shannon himself, but also Solomon Kullback (1907–1994), and Alan M. Turing (1912–1954) were deeply involved as cryptanalysts throughout WWII. Turing is best known for his role in theoretical computer science; his statistical thinking has been analysed in (Good 1979). Kullback's experience as a cryptanalyst lead him to write one of the first modern textbooks of statistical Cryptology (Kullback 1976).

**Conclusion: polymaths of past centuries**

The excessive specialization of modern science has induced an anachronistic tendency to associate each name to a discipline: Al Kindī to Philosophy, Alberti to Architecture, Viète to Algebra, Lieber to Linguistics, Babbage and Turing to Computer Science, Markov to Probability, just to mention a few names from the previous sections. The contents of this article would probably have puzzled most of its actors, who perceived knowledge as a whole and viewed themselves as scholars, or philosophers, in a very broad sense of the word: see (Burke 2012, chapter 6) on the progressive division of knowledge during the second half of the 19th century. The impressive variety of subjects treated by Al Kindī in his different treatises makes it meaningless to lock him into a single discipline (Mrayati et al. 2002, 78). The same can be said of Alberti (Grafton 2000), or Babbage

(Babbage 1864).

Most of the characters of our story had a training in classics, and knew Latin and Greek. Significantly, when Quetelet reviews Hayez's thesis in (Quetelet 1829, 391–394), he relays the new doctor's complaint, about having to write in Latin on a subject unheard of in Rome. Poetry, particularly in Latin was familiar to all. A well known anecdote reports that Leonhard Euler (1707–1783), was able to recite Virgil's Aeneid from any page he was given the number of. The statistician Francis Y. Edgeworth (1845–1926) made a statistical study of the dactyls in the same Aeneid, at the time when authorship attribution studies were just beginning (Stigler 1999, 111).

Cryptography, if only as a hobby, has also been familiar to most educated people for a very long time. Many well known names could be cited, like Bacon, Newton, Casanova, Euler, as being more or less marginally associated to cryptograms. "The Doctor" John Wallis (1616–1703), epitomizes our subject. He was an active cryptanalyst during the Glorious Revolution: he left examples of deciphered correspondence, together with instructions, translated into English in (Davys 1737). He was one of the very best mathematicians of his time; in Statistics he contributed the first English book on probabilities, shortly before the introduction of the first mortality tables by John Graunt and Edmund Halley. Last but not least, he wrote the first grammar of the English language… in Latin.

Eclectism, love of humanities and taste for riddles, we believe that these were important instigators for the mutual enrichment of Cryptology, Quantitative Analysis, and Statistics, through letter counting. There may have been, at least in some cases, another marginal factor. Most of the figures cited in the article were professors as well. Cryptograms make excellent classroom exercises. Letter counting also is a good collaborative activity to make students understand data collection (Richardson et al. 2004). We find it significant that Markov's second study on Aksakov only appeared in his textbook on probability (Link 2006a). Think of a time when some had to teach statistics without computers nor internet to download data sets from…